\documentclass[12pt,leqno]{article}

\usepackage[active]{srcltx}

\usepackage{amsfonts,latexsym,amsmath,amssymb,amsthm}
\usepackage{fullpage}
\newtheorem{theorem}{Theorem}[section]
\newtheorem{lemma}[theorem]{Lemma}

\newtheorem{proposition}[theorem]{Proposition}
\newtheorem{corollary}[theorem]{Corollary}
\newtheorem{remark}[theorem]{Remark}

\newtheorem{question}[theorem]{Question}
\theoremstyle{definition}

\newtheorem*{example}{Example}

\theoremstyle{remark}

\newtheorem*{note*}{Note}

\numberwithin{equation}{section}

\makeatletter
\def\R{{\mathbb R}}
\def\N{{\mathbb N}}
\def\e{\varepsilon}
\newcommand{\rank}{\mathop{\operator@font rank}}
\newcommand{\conv}{\mathop{\operator@font conv}}
\newcommand{\vol}{\mathop{\operator@font vol}}
\newcommand{\onetagright}{\tagsleft@false}
\makeatother

\newcommand{\ls}{\leqslant}
\newcommand{\gr}{\geqslant}
\renewcommand{\epsilon}{\varepsilon}
\newcommand{\prend}{$\quad \hfill \Box$}
\usepackage{color}

\usepackage{hyperref}

\begin{document}
\small

\title{\bf Variants of the Busemann-Petty problem and of the Shephard problem}

\medskip

\author{Apostolos Giannopoulos and Alexander Koldobsky}

\date{}

\maketitle

\begin{abstract}
\footnotesize We provide an affirmative answer to a variant of the Busemann-Petty problem, proposed by V.~Milman:
Let $K$ be a convex body in ${\mathbb R}^n$ and let $D$ be a compact subset of ${\mathbb R}^n$ such that,
for some $1\ls k\ls n-1$,
\begin{equation*}|P_F(K)|\ls |D\cap F|\end{equation*}
for all $F\in G_{n,k}$, where $P_F(K)$ is the orthogonal projection of $K$ onto $F$ and $D\cap F$ is the
intersection of $D$ with $F$. Then,
\begin{equation*}|K|\ls |D|.\end{equation*}
We also provide estimates for the lower dimensional Busemann-Petty and Shephard problems, and
we prove separation in the original Busemann-Petty problem.
\end{abstract}

\section{Introduction}

The Busemann-Petty problem was posed in \cite{Busemann-Petty-1956},
first in a list of ten problems concerning central sections of symmetric convex bodies in ${\mathbb R}^n$
and coming from questions in Minkowski geometry. It was originally formulated as
follows:
\begin{quote}{\sl Assume that $K$ and $D$ are origin-symmetric convex bodies in ${\mathbb R}^n$ and satisfy
\begin{equation}\label{eq:intro-1}|K\cap\xi^{\perp }|\ls |D\cap\xi^{\perp }|\end{equation}
for all $\xi\in S^{n-1}$. Does it follow that $|K|\ls |D|$?}
\end{quote}
Here $\xi^\perp$ is the central hyperplane perpendicular to $\xi $.
The answer is affirmative if $n\ls 4$ and negative if $n\gr 5$ (for the history and the solution to this problem,
see the monographs \cite{Gardner-book} and \cite{Koldobsky-book}). The isomorphic version of the Busemann-Petty problem
asks if there exists an absolute constant $C_1>0$ such that whenever $K$ and $D$ satisfy \eqref{eq:intro-1} we have $|K|\ls C_1|D|$.
This question is equivalent to the slicing problem and to the isotropic constant conjecture asking if
\begin{equation}\label{eq:intro-2}L_n:= \max\{ L_K:K\ \hbox{is isotropic in}\ {\mathbb R}^n\}\end{equation} is a bounded sequence.
More precisely, it is known that if $K$ and $D$ are two centered convex bodies
in ${\mathbb R}^n$ such that \eqref{eq:intro-1} holds true for all $\xi\in S^{n-1}$, then
\begin{equation}\label{eq:intro-3}|K|^{\frac{n-1}{n}}\ls c_1L_n\,|D|^{\frac{n-1}{n}},\end{equation}
where $c_1>0$ is an absolute constant. Regarding $L_n$, Bourgain proved in \cite{Bourgain-1991} that $L_n\ls
c\sqrt[4]{n}\log\! n$, and Klartag \cite{Klartag-2006} improved this bound to $L_n\ls c\sqrt[4]{n}$. A second proof of Klartag's bound
appears in \cite{Klartag-EMilman-2012}. For more information on isotropic convex bodies and log-concave measures see \cite{BGVV-book}.

Shephard's problem (see \cite{Shephard-1964}) is dual to the Busemann-Petty problem.
\begin{quote}{\sl Let $K$ and $D$ be two centrally symmetric convex bodies in ${\mathbb R}^n$. Suppose that
\begin{equation}\label{eq:intro-4}|P_{\xi^{\perp } }(K)|\ls |P_{\xi^{\perp} }(D)|\end{equation}
for every $\xi \in S^{n-1}$, where $P_{\xi^{\perp }}(A)$ is the orthogonal projection of $A\subset {\mathbb R}^n$
onto $\xi^{\perp }$. Does it follow that $|K|\ls |D|$?}
\end{quote}

The answer is affirmative if $n=2$, but shortly after it was posed, Shephard's question was answered in the negative for all $n\gr 3$.
This was done independently by Petty in \cite{Petty-1967} who gave an explicit counterexample in ${\mathbb R}^3$, and by Schneider
in \cite{Schneider-1967} for all $n\gr 3$. After these counterexamples, one might try to relax the question, asking for the smallest constant
$C_n$ (or the order of growth of this constant $C_n$ as $n\to \infty $) for which: if $K,D$ are centrally symmetric convex bodies in
${\mathbb R}^n$ and $|P_{\xi^{\perp }}(K)|\ls |P_{\xi^{\perp }}(D)|$ for all $\xi \in S^{n-1}$ then $|K|\ls C_n|D|$.

Such a constant $C_n$ does exist, and a simple argument, based on John's theorem, shows that $C_n\ls c\sqrt{n}$,
where $c>0$ is an absolute constant. On the other hand, K. Ball has proved in \cite{Ball-1991} that this simple
estimate is optimal: one has $C_n\simeq\sqrt{n}$.

\smallskip

In the first part of this note we discuss a variant of the two problems, proposed by V.~Milman at the Oberwolfach meeting on Convex Geometry
and its Applications (December 2015):

\begin{question}[V. Milman]\label{question:vitali}Assume that $K$ and $D$ are origin-symmetric convex bodies in ${\mathbb R}^n$ and satisfy
\begin{equation}\label{eq:intro-5}|P_{\xi^{\perp }}(K)|\ls |D\cap\xi^{\perp }|\end{equation}
for all $\xi\in S^{n-1}$. Does it follow that $|K|\ls |D|$?
\end{question}

In Section \ref{sect2} we show that the answer to this question is affirmative. In fact, the lower dimensional analogue of the problem has an affirmative
answer. Moreover, one can drop the symmetry assumptions and even the assumption of convexity for $D$.

\begin{theorem}\label{th:intro-1}Let $K$ be a convex body in ${\mathbb R}^n$ and let $D$ be a compact subset of ${\mathbb R}^n$ such that,
for some $1\ls k\ls n-1$,
\begin{equation}\label{eq:intro-6}|P_F(K)|\ls |D\cap F|\end{equation}
for all $F\in G_{n,n-k}$. Then,
\begin{equation}\label{eq:intro-7}|K|\ls |D|.\end{equation}
\end{theorem}
We also prove stability and separation in Theorem \ref{th:intro-1}.
In the hyperplane case, and assuming that $K$ and $D$ are centered convex bodies, i.e. their center of mass is at the origin,
we can provide a more precise answer in terms of the isotropic constant $L_D$ of $D$.

\begin{theorem}\label{th:intro-2}Let $K$ and $D$ be two centered convex bodies in ${\mathbb R}^n$ such that
\begin{equation}\label{eq:intro-8}|P_{\xi^{\perp }}(K)|\ls |D\cap\xi^{\perp }|\end{equation}
for all $\xi\in S^{n-1}$. Then,
\begin{equation}\label{eq:intro-9}|K|\ls \frac{c}{L_D}\,|D|,\end{equation}
where $c>0$ is an absolute constant. \end{theorem}

This means that if the hyperplane conjecture is not true then one can even have ``pathologically good" (with respect to Question \ref{question:vitali})
pairs of convex bodies. The proof of Theorem \ref{th:intro-2} carries over to higher codimensions but the dependence on $L_D$ becomes more
complicated and we prefer not to include the full statement of this version.

\bigbreak

In Section \ref{sect3} we collect some estimates on the lower dimensional Busemann-Petty problem. Let $1\ls k\ls n-1$
and let $\beta_{n,k}$ be the smallest constant $\beta >0$ with the following property: For
every pair of centered convex bodies $K$ and $D$ in ${\mathbb R}^n$ that satisfy
\begin{equation}\label{eq:intro-10}|K\cap F|\ls |D\cap F|\end{equation}
for all $F\in G_{n,n-k}$, one has
\begin{equation}\label{eq:intro-11}|K|^{\frac{n-k}{n}}\ls \beta^k\,|D|^{\frac{n-k}{n}}.\end{equation}
The following question is open:

\begin{question}\label{question:low-dim-BP}{\sl Is it true that there exists an absolute constant $C_2>0$ such that $\beta_{n,k}\ls C_2$ for all $n$ and $k$?}
\end{question}

Bourgain and Zhang \cite{BZ} showed that $\beta_{n,k}>1$ if $n-k>3.$ It is not known whether $\beta_{n,k}$
has to be greater than 1 when $n\ge 5$ and $n-k=2$ or $n-k=3.$
It was proved in \cite{Koldobsky-4} and by a different method in \cite{Chasapis-Giannopoulos-Liakopoulos-2015}
that $\beta_{n,k}\le C\sqrt{n/k}(\log(en/k))^{3/2},$ where $C$ is an absolute constant. In this note, we observe that the answer to Question \ref{question:low-dim-BP} is affirmative if the convex body $K$ has bounded isotropic constant, as follows.

\begin{theorem}\label{th:intro-3}Let $1\ls k\ls n-1$ and let $K$ be a centered convex body in ${\mathbb R}^n$ and $D$ a compact subset
of ${\mathbb R}^n$ such that
\begin{equation}\label{eq:intro-12}|K\cap F|\ls |D\cap F|\end{equation}
for all $F\in G_{n,n-k}$. Then,
\begin{equation}\label{eq:intro-13}|K|^{\frac{n-k}{n}}\ls (c_0L_K)^k\,|D|^{\frac{n-k}{n}}.\end{equation}
where $c_0>0$ is an absolute constant. \end{theorem}

Theorem \ref{th:intro-3} is a refinement of the estimate $\beta_{n,k}\ls cL_n$, which was shown in \cite{Chasapis-Giannopoulos-Liakopoulos-2015}.
The proof is based on estimates from \cite{Dafnis-Paouris-2012} and on Grinberg's inequality (see \eqref{eq:main-3} in Section 2).

We also discuss the lower dimensional Shephard problem. Let $1\ls k\ls n-1$
and let $S_{n,k}$ be the smallest constant $S >0$ with the following property: For
every pair of convex bodies $K$ and $D$ in ${\mathbb R}^n$ that satisfy
\begin{equation}\label{eq:intro-S-1}|P_F(K)|\ls |P_F(D)|\end{equation}
for all $F\in G_{n,n-k}$, one has
\begin{equation}\label{eq:intro-S-2}|K|^{\frac{1}{n}}\ls S\,|D|^{\frac{1}{n}}.\end{equation}

\begin{question}\label{question:low-dim-S}{\sl Is it true that there exists an absolute constant $C_3>0$ such that $S_{n,k}\ls C_3$ for all $n$ and $k$?}
\end{question}
Goodey and Zhang \cite{GZ} proved that $S_{n,k}>1$ if $n-k>1.$
In Section \ref{sect4} we prove the following result.

\begin{theorem}\label{th:intro-S}Let $K$ and $D$ be two convex bodies in ${\mathbb R}^n$ such that
\begin{equation}|P_F(K)|\ls |P_F(D)|\end{equation}
for every $F\in G_{n,n-k}$. Then,
\begin{equation}|K|^{\frac{1}{n}}\ls c_1\sqrt{\frac{n}{n-k}}\log\left (\frac{en}{n-k}\right )\,|D|^{\frac{1}{n}},\end{equation}
where $c_1>0$ is an absolute constant. It follows that $S_{n,k}$ is bounded by an absolute constant if $\frac{k}{n-k}$ is bounded.
\end{theorem}

We also prove a general estimate, which is logarithmic in $n$ and valid for all $k$. The proof is based on estimates from \cite{Paouris-Pivovarov-2013}.

\begin{theorem}Let $K$ and $D$ be two convex bodies in ${\mathbb R}^n$ such that
\begin{equation}|P_F(K)|\ls |P_F(D)|\end{equation}
for every $F\in G_{n,n-k}$. Then,
\begin{equation}|K|^{\frac{1}{n}}\ls \frac{c_1\,\min w(\tilde{D})}{\sqrt{n}}\,|D|^{\frac{1}{n}}\ls c_2(\log n)|D|^{\frac{1}{n}},\end{equation}
where $c_1,c_2>0$ are absolute constants, $w(A)$ is the mean width of a centered convex body $A$, and the minimum is over all linear
images $\tilde{D}$ of $D$ that have volume $1$.
\end{theorem}

Lutwak \cite{Lutwak-1988} proved that the answer to the Busemann-Petty problem is affirmative if the body $K$ with smaller sections belongs to a special class
of intersection bodies; see definition below. In Section \ref{sect5} we prove separation in the Busemann-Petty
problem, which can be considered as a refinement of Lutwak's result.

\begin{theorem} \label{main-int} Suppose that $\e>0$,  $K$ and $D$ are origin-symmetric
star bodies in $\R^n,$ $K$ is an intersection body. If
\begin{eqnarray}\label{sect1}
|K\cap \xi^\bot| \ls  |D\cap \xi^\bot| - \e,
\end{eqnarray}
for every $\xi\in S^{n-1}$, then
$$|K|^{\frac{n-1}n}  \ls  |D|^{\frac{n-1}n} - c\e  \frac{1}{\sqrt{n}M(\overline{K})},$$
where $c>0$ is an absolute constant and $\overline{K}=|K|^{-\frac{1}{n}}K$.
\end{theorem}

Note that if $\overline{K}$ is convex isotropic then \begin{equation}\frac{1}{M(\overline{K})} \gr\ c_1 \frac{n^{1/10}L_K}{\log^{2/5}(e+n)}\gr c_2\frac{n^{1/10}}{\log^{2/5}(e+n)}\end{equation}
and if $\overline{K}$ is convex and is in the minimal mean width position then we have
\begin{equation}\frac{1}{M(\overline{K})} \gr\ c_3\frac{\sqrt{n}}{\log (e+n)},\end{equation}
so the constant in Theorem \ref{main-int} does not depend
on the bodies. This is an improvement of a previously known result from \cite{K??}. Also note that
stability in the Busemann-Petty problem is easier and was proved in \cite{K??}, as follows.
If $K$ is an intersection body in $\R^n,$ $D$ is an origin-symmetric star
body in $\R^n$ and $\e>0$ so that
\begin{equation}|K\cap \xi^\bot|\ls |D\cap \xi^\bot|+\e \end{equation}
for every $\xi\in S^{n-1}$, then
\begin{equation}|K|^{\frac{n-1}n}\ls |L|^{\frac{n-1}n} + c_n\e,\end{equation}
where $c_n=|B_2^{n-1}|/|B_2^n|^{\frac{n-1}n} < 1.$
The constant is optimal. For more results on stability and separation in volume comparison
problems and for applications of such results, see \cite{Koldobsky-3}.

\section{Milman's variant of the two problems}\label{sect2}

We work in ${\mathbb R}^n$, which is equipped with a Euclidean structure $\langle\cdot ,\cdot\rangle $. We denote the corresponding
Euclidean norm by $\|\cdot \|_2$, and write $B_2^n$ for the Euclidean unit ball, and $S^{n-1}$ for the unit sphere. Volume is
denoted by $|\cdot |$. We write $\omega_n$ for the volume of $B_2^n$ and $\sigma $ for the rotationally invariant probability measure on
$S^{n-1}$. We also denote the Haar measure on $O(n)$ by $\nu $. The Grassmann manifold $G_{n,m}$ of $m$-dimensional subspaces of
${\mathbb R}^n$ is equipped with the Haar probability measure $\nu_{n,m}$. Let $1\ls m\ls n-1$ and $F\in G_{n,m}$. We will denote the
orthogonal projection from $\mathbb R^{n}$ onto $F$ by $P_F$. We also define $B_F=B_2^n\cap F$ and $S_F=S^{n-1}\cap F$.

The letters $c,c^{\prime }, c_1, c_2$ etc. denote absolute positive constants whose value may change from line to line. Whenever we
write $a\simeq b$, we mean that there exist absolute constants $c_1,c_2>0$ such that $c_1a\ls b\ls c_2a$.  Also if $K,L\subseteq
\mathbb R^n$ we will write $K\simeq L$ if there exist absolute constants $c_1, c_2>0$ such that $c_{1}K\subseteq L \subseteq
c_{2}K$.

A convex body in ${\mathbb R}^n$ is a compact convex subset $K$ of
${\mathbb R}^n$ with nonempty interior. We say that $K$ is origin-symmetric if $K=-K$. We say that $K$ is centered if
the center of mass of $K$ is at the origin, i.e.~$\int_K\langle
x,\theta\rangle \,d x=0$ for every $\theta\in S^{n-1}$. We denote by ${\cal K}_n$ the class of centered
convex bodies in ${\mathbb R}^n$. The support
function of $K$ is defined by $h_K(y):=\max \bigl\{\langle x,y\rangle :x\in K\bigr\}$, and
the mean width of $K$ is the average
\begin{equation}\label{eq:not-2}w(K):=\int_{S^{n-1}}h_K(\theta )\,d\sigma (\theta )\end{equation}
of $h_K$ on $S^{n-1}$. For basic facts from the Brunn-Minkowski theory and
the asymptotic theory of convex bodies we refer to the books \cite{Schneider-book} and \cite{AGA-book} respectively.

\smallskip

The proof of Theorem \ref{th:intro-1} is based on two classical results:
\begin{enumerate}
\item {\sl Aleksandrov's inequalities.} If $K$ is a convex body in ${\mathbb R}^n$ then the sequence
\begin{equation}\label{eq:main-1}Q_k(K)=\left
(\frac{1}{\omega_k}\int_{G_{n,k}}|P_F(K)|\,d\nu_{n,k}(F)\right
)^{1/k}\end{equation}is decreasing in $k$. This is a consequence of the Aleksandrov-Fenchel inequality (see \cite{Burago-Zalgaller-book}
and \cite{Schneider-book}). In particular, for every $1\ls k\ls n-1$ we have
\begin{equation}\label{eq:main-2}\left (\frac{|K|}{\omega_n}\right )^{\frac{1}{n}}\ls \left (\frac{1}{\omega_k}\int_{G_{n,k}}|P_F(K)|\,d\nu_{n,k}(F)\right )^{\frac{1}{k}}\ls w(K),\end{equation}
where $w(K)$ is the mean width of $K$.
\item {\sl Grinberg's inequality.}  If $D$ is a compact set in ${\mathbb R}^n$ then, for any $1\ls k\ls n-1$,
\begin{equation}\label{eq:main-3}\tilde{R}_k(D):=\frac{1}{|D|^{n-k}}\int_{G_{n,n-k}}|D\cap F|^n\,d\nu_{n,n-k}(F)\ls \frac{1}{|B_2^n|^{n-k}}\int_{G_{n,n-k}}|B_2^n\cap F|^n\,d\nu_{n,n-k}(F),\end{equation}
where $B_2^m$ is the Euclidean ball in ${\mathbb R}^m$ and $\omega_m=|B_2^m|$. This fact was proved by Grinberg in \cite{Grinberg-1990}.
It is useful to note that
\begin{equation}\label{eq:main-4}\tilde{R}_k(B_2^n):=\frac{\omega_{n-k}^n}{\omega_n^{n-k}}\ls e^{\frac{kn}{2}}.\end{equation}
Moreover, Grinberg proved that the quantity $\tilde{R}_k(D)$ on the left hand side of \eqref{eq:main-3} is invariant under $T\in GL(n)$: one has
\begin{equation}\label{eq:main-5}\tilde{R}_k(T(D))=\tilde{R}_k(D)\end{equation}
for every $T\in GL(n)$.
\end{enumerate}

\smallskip

\noindent {\bf Proof of Theorem \ref{th:intro-1}.} Let $K$ be a convex body in ${\mathbb R}^n$ and $D$ be a compact subset of ${\mathbb R}^n$. Assume that for some
$1\ls k\ls n-1$ we have
\begin{equation}\label{eq:main-6}|P_F(K)|\ls |D\cap F|\end{equation}
for all $F\in G_{n,n-k}$. From \eqref{eq:main-2} we get
\begin{equation}\label{eq:main-7}\left (\frac{|K|}{\omega_n}\right )^{\frac{n-k}{n}}\ls \frac{1}{\omega_{n-k}}\int_{G_{n,n-k}}|P_F(K)|\,d\nu_{n,n-k}(F).\end{equation}
Our assumption, H\"{o}lder's inequality and Grinberg's inequality give
\begin{align}\label{eq:main-8}&\frac{1}{\omega_{n-k}}\int_{G_{n,n-k}}|P_F(K)|\,d\nu_{n,n-k}(F) \ls \frac{1}{\omega_{n-k}}\int_{G_{n,n-k}}|D\cap F|\,d\nu_{n,n-k}(F)\\
\nonumber &\hspace*{1.5cm} \ls \frac{1}{\omega_{n-k}}\left (\int_{G_{n,n-k}}|D\cap F|^n\,d\nu_{n,n-k}(F)\right )^{\frac{1}{n}}\\
\nonumber &\hspace*{1.5cm}\ls  \frac{1}{\omega_{n-k}}\,\frac{\omega_{n-k}}{\omega_n^{\frac{n-k}{n}}}|D|^{\frac{n-k}{n}}=  \left (\frac{|D|}{\omega_n}\right )^{\frac{n-k}{n}}.
\end{align}
Therefore, $|K|\ls |D|$. \prend

\medskip

\begin{remark}\label{stability}\rm Slightly modifying the proof of Theorem \ref{th:intro-1} one can get
stability and separation results, as follows. Let $\varepsilon>0,$ and let $K$ and $D$ be as in Theorem \ref{th:intro-1}.
Suppose that for every $F\in G_{n,n-k}$
$$|P_F(K)| \le |D\cap F| \pm \varepsilon.$$
Then
$$|K|^{\frac{n-k}n} \le |D|^{\frac{n-k}n} \pm \gamma_{n,k}\varepsilon,$$
where
$\gamma_{n,k} = \frac{\omega_n^{\frac{n-k}{n}}}{\omega_{n-k}}\in (e^{-k/2}, 1).$ The plus sign corresponds to stability,
minus - to separation. Assuming that $\varepsilon =\max_F (|P_F(K) - |D\cap F|)$ in the stability result, we get
$$|K|^{\frac{n-k}n} - |D|^{\frac{n-k}n} \le \gamma_{n,k} \max_F (|P_F(K) - |D\cap F|).$$
On the other hand, if $\varepsilon =\min_F (|D\cap F| - |P_F(K)|)$ in the separation result, then
$$|D|^{\frac{n-k}n} - |K|^{\frac{n-k}n} \ge \gamma_{n,k} \min_F (|D\cap F| - |P_F(K)|).$$
\end{remark}

\medskip


\section{Estimates for the lower dimensional Busemann-Petty problem}\label{sect3}

In this section we provide some estimates for the lower dimensional Busemann-Petty problem. We need the next lemma, in which we collect known estimates about the quantities
\begin{equation}G_{n,k}(A):=\left (\int_{G_{n,n-k}}|A\cap F|^n\,d\nu_{n,n-k}(F)\right )^{\frac{1}{kn}},\end{equation}
where $A$ is a centered convex body in ${\mathbb R}^n$. The proofs of \eqref{eq:main-9} and \eqref{eq:main-10} can be found
in \cite{Dafnis-Paouris-2012}, while \eqref{eq:main-11} follows from \eqref{eq:main-3} and \eqref{eq:main-4}.

\begin{lemma}\label{lem:main-2}Let $A$ be a centered convex body in ${\mathbb R}^n$. Then,
\begin{equation}\label{eq:main-9}\frac{c_1}{L_A}|A|^{\frac{n-k}{kn}}\ls G_{n,k}(A)\ls \frac{c_2L_k}{L_A}|A|^{\frac{n-k}{kn}}\ls \frac{c_3\sqrt[4]{k}}{L_A}|A|^{\frac{n-k}{kn}}.\end{equation}
Moreover,
\begin{equation}\label{eq:main-10}G_{n,k}(A)\ls c_4\sqrt{n/k}\,(\log (en/k))^{\frac{3}{2}}|A|^{\frac{n-k}{kn}}.\end{equation}
Finally, for every compact subset $D$ of ${\mathbb R}^n$ we have
\begin{equation}\label{eq:main-11}G_{n,k}(D)\ls\sqrt{e}|D|^{\frac{n-k}{kn}}.\end{equation}
\end{lemma}

Using Lemma \ref{lem:main-2} we show that the lower dimensional Busemann-Petty problem (Question \ref{question:low-dim-BP})
has an affirmative answer if the body $K$ has bounded isotropic constant.

\medskip

\noindent {\bf Proof of Theorem \ref{th:intro-3}.} Since $|K\cap F|\ls |D\cap F|$ for all $F\in G_{n,n-k}$, we know that
\begin{equation}\label{eq:LBP-1}G_{n,k}(K)\ls G_{n,k}(D).\end{equation}
Using \eqref{eq:main-9} and \eqref{eq:main-11} we write
\begin{equation}\label{eq:LBP-2}\frac{c_1}{L_K}|K|^{\frac{n-k}{kn}} \ls G_{n,k}(K)\ls G_{n,k}(D)
\ls\sqrt{e}|D|^{\frac{n-k}{kn}},\end{equation}
and the result follows. \prend.

\begin{remark}\label{rem:LBP-1}\rm Theorem \ref{th:intro-3} shows that if $K$ belongs to the class
\begin{equation}\label{eq:class-alpha}{\cal K}_n(\alpha ):=\{K\in {\cal K}_n: L_K\ls\alpha\}\end{equation}
for some $\alpha >0$, then for every compact set $D$ in ${\mathbb R}^n$ which satisfies $|K\cap F|\ls |D\cap F|$ for all $F\in G_{n,n-k}$ we have
\begin{equation}\label{eq:LBP-5}|K|^{\frac{n-k}{n}}\ls (c_0\alpha )^k\,|D|^{\frac{n-k}{n}}.\end{equation}
Classes of convex bodies with uniformly bounded isotropic constant include: unconditional convex bodies,
convex bodies whose polar bodies contain large affine cubes, the
unit balls of $2$-convex spaces with a given constant
$\alpha $, bodies with small diameter (in particular, the class of
zonoids) and the unit balls of the Schatten classes (see \cite[Chapter 4]{BGVV-book}).\end{remark}

\begin{example}\label{ex:beta-lower}\rm K. Ball has proved in \cite{Ball-1989} that for every $1\ls k\ls n-1$ and every $F\in G_{n,n-k}$
we have
\begin{equation}|Q_n\cap F|\ls 2^{\frac{k}{2}},\end{equation}
where $Q_n$ is the cube of volume $1$ in ${\mathbb R}^n$. Consider the ball $B_{n,k}=r_{n,k}B_2^n$, where
\begin{equation}\omega_{n-k}r_{n,k}^{n-k}=2^{\frac{k}{2}}.\end{equation}
Then, for every $F\in G_{n,n-k}$ we have
\begin{equation}|Q_n\cap F|\ls |B_{n,k}\cap F|.\end{equation}
Therefore,
\begin{equation}1=|Q_n|\ls \beta_{n,k}^k|B_{n,k}|^{\frac{n-k}{n}}=\beta_{n,k}^k\omega_n^{\frac{n-k}{n}}r_{n,k}^{n-k}=2^{\frac{k}{2}}\beta_{n,k}^k\frac{\omega_n^{\frac{n-k}{n}}}{\omega_{n-k}}.\end{equation}
This proves that
\begin{equation}\beta_{n,k}\gr \frac{1}{\sqrt{2}}\left (\frac{\omega_{n-k}}{\omega_n^{\frac{n-k}{n}}}\right )^{\frac{1}{k}}\sim \frac{1}{\sqrt{2}}\left (\frac{n}{n-k}\right )^{\frac{n-k+1}{2k}}\end{equation}
as $n,k\to\infty $. Fix $d\gr 2$ and consider $n$ and $k$ that satisfy $n=(d+1)k$. Then, we have the following:
\end{example}

\begin{proposition}\label{prop:beta-lower-bound}For every $d\gr 2$ there exists $k(d)\in {\mathbb N}$ such that
\begin{equation}\beta_{(d+1)k,k}\gr \frac{1}{\sqrt{2}}\left (1+\frac{1}{d}\right )^{\frac{d}{2}}>1\end{equation}
for all $k\gr k(d)$. \prend \end{proposition}

\medskip

A variant of the proof of Theorem \ref{th:intro-3} (based again on Lemma \ref{lem:main-2}) establishes Theorem \ref{th:intro-2}.

\medskip

\noindent {\bf Proof of Theorem \ref{th:intro-2}.} Let $K$ be a convex body in ${\mathbb R}^n$ and $D$ be a compact subset of ${\mathbb R}^n$
such that $|P_{\xi^{\perp }}(K)|\ls |D\cap\xi^{\perp }|$ for every $\xi\in S^{n-1}$. From Lemma \ref{lem:main-2} we know that
\begin{equation}\label{eq:main-12}G_{n,1}(D)\ls \frac{c_1}{L_D}|D|^{\frac{n-1}{n}},\end{equation}
where $c_1>0$ is an absolute constant. Then,
\begin{align}\label{eq:main-13}
\left (\frac{|K|}{\omega_n}\right )^{\frac{n-1}{n}} &\ls \frac{1}{\omega_{n-1}}\int_{S^{n-1}}|P_{\xi^{\perp }}(K)|\,d\sigma (\xi )
\ls \frac{1}{\omega_{n-1}}\int_{S^{n-1}}|D\cap \xi^{\perp }|\,d\sigma (\xi )\\
\nonumber &= \frac{1}{\omega_{n-1}}\left (\int_{S^{n-1}}|D\cap \xi^{\perp }|^n\,d\sigma (\xi )\right )^{\frac{1}{n}}\\
\nonumber &= \frac{1}{\omega_{n-1}}G_{n,1}(D)\ls \frac{c_1}{\omega_{n-1}L_D}|D|^{\frac{n-1}{n}},
\end{align}
which implies that
\begin{equation}\label{eq:main-14}|K|\ls \frac{c_2\omega_n}{(\omega_{n-1}L_D)^{\frac{n}{n-1}}}\,|D|\ls \frac{c_3}{L_D}\,|D|,\end{equation}
where $c_2,c_3>0$ are absolute constants. \prend

\section{Estimates for the lower dimensional Shephard problem}\label{sect4}

In this section we discuss the lower dimensional Shephard problem. First, we recall some facts for
the class of zonoids. A zonoid is a limit of Minkowski sums of line segments in the Hausdorff metric. Equivalently, a
symmetric convex body $Z$ is a zonoid if and only if its polar body is the unit ball of an $n$-dimensional subspace of an $L_1$-space;
i.e. if there exists a positive measure $\mu $ (the supporting measure of $Z$) on $S^{n-1}$ such that
\begin{equation*}h_Z(x)=\| x\|_{Z^{\circ }}=\frac{1}{2}\int_{S^{n-1}}|\langle x,y\rangle |d\mu (y).\end{equation*}
The class of origin-symmetric zonoids coincides with the class of projection bodies. Recall that the projection body $\Pi K$ of a convex body $K$ is the
symmetric convex body whose support function is defined by
\begin{equation*}h_{\Pi K} (\xi )=|P_{\xi^{\perp } }(K)|, \qquad \xi\in S^{n-1}.\end{equation*}
From Cauchy's formula
\begin{equation*}|P_{\xi^{\perp } }(K)|={\frac{1}{2}}\int_{S^{n-1}}|\langle u,\xi \rangle |\;d\sigma_K(u),\end{equation*}
where $\sigma_K$ is the surface area measure of $K$, it follows that the projection body of $K$ is a zonoid whose
supporting measure is $\sigma_K$. Minkowski's existence theorem implies that, conversely, every zonoid
is the projection body of some symmetric convex body in ${\mathbb R}^n$.

Zonoids play a central role in the study of the original Shephard problem: suppose that $K$ is a convex body in ${\mathbb R}^n$
and $Z$ is a zonoid in ${\mathbb R}^n$, and that $|P_{\xi^{\perp }}(K)|\ls |P_{\xi^{\perp }}(Z)|$ for all $\xi \in S^{n-1}$. Then,
\begin{equation}|K|\ls |Z|.\end{equation}
The proof involves writing $Z=\Pi D$ for some convex body $D$, using the identity $V_1(K,\Pi D)=V_1(D,\Pi K)$ (where $V_1(A,B)$ is the
mixed volume $V(A,\ldots ,A,B)$), the hypothesis in the form $\Pi (K)\subseteq \Pi (Z)$, and the
monotonicity of $V_1(D,.)$, to write
\begin{equation*}|Z|=V_1(Z,Z)=V_1(Z,\Pi D)=V_1(D,\Pi Z)\gr V_1(D,\Pi (K))=
V_1(K,\Pi D)=V_1(K,Z)\gr |K|^{{\frac{n-1}{n}}}|Z|^{\frac{1}{n}},\end{equation*}
where in the last step we also employ Minokowski's first inequality. This shows that $|Z|\gr |K|$.

Since any projection of a zonoid is a zonoid, using an inductive argument we can prove the following (for a detailed account
on this topic, see \cite[Chapter 4]{Gardner-book}).

\begin{theorem}\label{th:SZ-1}Let $K$ be a convex body and let $Z$ be a zonoid in ${\mathbb R}^n$ such that
\begin{equation}|P_F(K)|\ls |P_F(Z)|\end{equation}
for every $F\in G_{n,n-k}$. Then,
\begin{equation}\label{eq:SZ-1}|K|\ls |Z|.\end{equation}
\end{theorem}

Using Theorem \ref{th:SZ-1} and the fact that every ellipsoid is a zonoid, we can give a simple bound for the constants $S_{n,k}$.

\begin{proposition}\label{prop:SZ-simple}For all $n$ and $1\ls k\ls n-1$ we have $S_{n,k}\ls c_0\sqrt{n}$, where $c_0>0$ is an absolute constant.\end{proposition}

\noindent {\it Proof.} Let $K$ and $D$ be two convex bodies in ${\mathbb R}^n$ such that $|P_F(K)|\ls |P_F(D)|$ for every $F\in G_{n,n-k}$. There
exists an ellipsoid ${\cal E}$ in ${\mathbb R}^n$ such that $D\subseteq {\cal E}$ and $|{\cal E}|^{\frac{1}{n}}\ls c_0\sqrt{n}\,|D|^{\frac{1}{n}}$,
where $c_0>0$ is an absolute constant (for example, see \cite{Ball-1991c} where a sharp estimate for $c_0$ is also given). Since $D\subseteq {\cal E}$,
we have
\begin{equation}|P_F(K)|\ls |P_F(D)|\ls |P_F({\cal E})|\end{equation}
for all $F\in G_{n,n-k}$. Since ${\cal E}$ is a zonoid, Theorem \ref{th:SZ-1} implies that
\begin{equation}|K|^{\frac{1}{n}}\ls |{\cal E}|^{\frac{1}{n}}\ls c_0\sqrt{n}\,|D|^{\frac{1}{n}}.\end{equation}
This shows that $S_{n,k}\ls c_0\sqrt{n}$. \prend

\medskip

We can elaborate on this argument if we use Pisier's theorem from \cite{Pisier-1989} on the existence of $\alpha $-regular $M$-ellipsoids for
symmetric convex bodies in ${\mathbb R}^n$ (see \cite[Theorem 1.13.3]{BGVV-book}):

\begin{theorem}[Pisier]\label{th:general-pisier-alpha-regular}
For every $0<\alpha <2$ and every symmetric convex body $A$ in ${\mathbb
R}^n$, there exists an ellipsoid ${\cal{E}_{\alpha }}$ such that
\begin{equation*}\max\{ N(A,t{\cal{E}_{\alpha }}),N({\cal{E}_{\alpha }},tA)\}\ls\exp
\left (\frac{c(\alpha )n}{t^{\alpha }}\right )\end{equation*} for
every $t\gr 1$, where $c(\alpha )$ is a constant depending only on
$\alpha $ and satisfies $c(\alpha )=O\big ((2-\alpha )^{-\alpha/2}\big )$ as
$\alpha\to 2$.
\end{theorem}

\begin{theorem}\label{th:SZ-3}Let $1\ls m\ls n-1$ and let $K$ and $D$ be two convex bodies in ${\mathbb R}^n$ such that
\begin{equation}|P_F(K)|\ls |P_F(D)|\end{equation}
for every $F\in G_{n,m}$. Then,
\begin{equation}|K|^{\frac{1}{n}}\ls c_1\sqrt{\frac{n}{m}}\log\left(\frac{en}{m}\right )\,|D|^{\frac{1}{n}},\end{equation}
where $c_1>0$ is an absolute constant.
\end{theorem}

\noindent {\it Proof.} Consider the difference body $D-D$ of $D$, and the ellipsoid ${\cal E}_{\alpha }$ from Theorem \ref{th:general-pisier-alpha-regular}, where $\alpha\in (0,2)$
will be chosen in the end, that corresponds to $A=D-D$. Note that
\begin{equation}N({\cal E}_{\alpha },c(\alpha )^{1/\alpha }(D-D))\ls e^n,\end{equation}
therefore
\begin{equation}\label{eq:basic}|{\cal E}_{\alpha }|^{\frac{1}{n}}\ls ec(\alpha )^{1/\alpha }|D-D|^{\frac{1}{n}}.\end{equation}
Since
\begin{equation}N(P_F(D-D), P_F(t{\cal E}_{\alpha }))\ls N(D-D,t{\cal{E}_{\alpha }})\ls \exp
\left (\frac{c(\alpha )n}{t^{\alpha }}\right )\end{equation}
for every $F\in G_{n,m}$, we have
\begin{equation}|P_F(D-D)|\ls \exp
\left (\frac{c(\alpha )n}{t_{n,m,\alpha }^{\alpha }}\right )|P_F(t_{n,m,\alpha }{\cal E}_{\alpha })|=|P_F(et_{n,m,\alpha }{\cal E}_{\alpha })|\end{equation}
if we choose
\begin{equation}t_{n,m,\alpha }=\left (\frac{c(\alpha )n}{m}\right )^{\frac{1}{\alpha }}.\end{equation}
Now, if we set ${\cal E}:=et_{n,m,\alpha }{\cal E}_{\alpha }$, we have
\begin{equation}|P_F(K)|\ls |P_F(D)|\ls |P_F(D-D)|\ls |P_F({\cal E})|\end{equation}
for every $F\in G_{n,m}$, and since ${\cal E}$ is a zonoid, Theorem \ref{th:SZ-1} shows that $|K|\ls |{\cal E}|$.
Using also \eqref{eq:basic} and the fact that $c(\alpha )=O\big ((2-\alpha )^{-\alpha/2}\big )$, we get
\begin{equation}|K|^{\frac{1}{n}}\ls et_{n,m,\alpha }|{\cal E}_{\alpha }|^{\frac{1}{n}}\ls e^2t_{n,m,\alpha }c(\alpha )^{1/\alpha }|D-D|^{\frac{1}{n}}
\ls \frac{c_1}{2-\alpha }\left (\frac{n}{m}\right )^{\frac{1}{\alpha }}|D|^{\frac{1}{n}}, \end{equation}
where $c_1>0$ is an absolute constant (we have also used the fact that $|D-D|^{\frac{1}{n}}\ls 4|D|^{\frac{1}{n}}$ by the Rogers-Shephard
inequality). Choosing $\alpha =2-\frac{1}{\log\left(\frac{en}{m}\right )}$ we get the result. \prend

\begin{remark}\label{rem:lutwak-conj}\rm The lower dimensional Shephard problem is related to Lutwak's conjectures about
the affine quermassintegrals: for every convex body $K$ in ${\mathbb R}^n$ and every $1\ls m\ls
n-1$, the quantities
\begin{equation}\Phi_{n-m}(K)=\frac{\omega_n}{\omega_m}\left
(\int_{G_{n,m}}|P_F(K)|^{-n}d\nu_{n,m}(F)\right )^{-1/n},\end{equation}
were introduced by Lutwak in \cite{Lutwak-1984} (and Grinberg proved in \cite{Grinberg-1990}
that these quantities are invariant under volume preserving affine transformations).
Lutwak conjectured in \cite{Lutwak-1988b} that the affine quermassintegrals satisfy the inequalities
\begin{equation}\omega_n^j\Phi_i(K)^{n-j}\ls \omega_n^i\Phi_j(K)^{n-i}\end{equation}
for all $0\ls i<j<n$, where we agree that $\Phi_0(K)=|K|$ and $\Phi_n(K)=\omega_n$.
Most of the conjectures about the affine quermassintegrals remain open
(see \cite[Chapter 9]{Gardner-book} for more details and references). If true, they
would imply the following (see also \cite{Dafnis-Paouris-2012}): there exist
absolute constants $c_1,c_2>0$ such that for every convex body $K$
in $\mathbb R^n$ and every $1\ls m \ls n-1$,
\begin{equation}\label{eq:phi-conj}c_1\sqrt{n/m}\,|K|^{\frac{1}{n}}\ls  \left (\int_{G_{n,m}}|P_F(K)|^{-n}\,d\nu_{n,m}(F)\right )^{-\frac{1}{mn}}\ls c_2\sqrt{n/m}\,|K|^{\frac{1}{n}}.\end{equation}
Assuming \eqref{eq:phi-conj} we can give an affirmative answer to Question \ref{question:low-dim-S}. Indeed, let $K$ and $D$ be two convex bodies in ${\mathbb R}^n$ such that
$|P_F(K)|\ls |P_F(D)|$ for every $F\in G_{n,n-k}$. We write
\begin{align}c_1\sqrt{n/(n-k)}\,|K|^{\frac{1}{n}}&\ls  \left (\int_{G_{n,n-k}}|P_F(K)|^{-n}\,d\nu_{n,n-k}(F)\right )^{-\frac{1}{(n-k)n}}\\
\nonumber &\ls \left (\int_{G_{n,n-k}}|P_F(D)|^{-n}\,d\nu_{n,n-k}(F)\right )^{-\frac{1}{(n-k)n}}\ls c_2\sqrt{n/(n-k)}\,|D|^{\frac{1}{n}},\end{align}
and this shows that $|K|^{\frac{1}{n}}\ls (c_2/c_1)\,|D|^{\frac{1}{n}}$.
\end{remark}

The left hand side of \eqref{eq:phi-conj} was proved by Paouris and Pivovarov in \cite{Paouris-Pivovarov-2013}:

\begin{theorem}[Paouris-Pivovarov]\label{th:main-3}Let $A$ be a convex body in ${\mathbb R}^n$. Then,
\begin{equation}\label{eq:main-17}\left (\int_{G_{n,m}}|P_F(A)|^{-n}\,d\nu_{n,m}(F)\right )^{-\frac{1}{mn}}\gr c\sqrt{n/m}|A|^{\frac{1}{n}}.\end{equation}
\end{theorem}

Using this fact one can obtain the following.

\begin{proposition}\label{prop:SP-1}Let $1\ls m\ls n-1$ and let $K$ and $D$ be two convex bodies in ${\mathbb R}^n$ such that
\begin{equation}|P_F(K)|\ls |P_F(D)|\end{equation} for every $F\in G_{n,m}$. Then,
\begin{equation}\label{eq:SP-1}|K|^{\frac{1}{n}}\ls \frac{c\,\min w(\tilde{D})}{\sqrt{n}}\,|D|^{\frac{1}{n}},\end{equation}
where $c>0$ is an absolute constant, $w(A)$ is the mean width of a centered convex body $A$, and the minimum is over all linear
images $\tilde{D}$ of $D$ that have volume $1$.
\end{proposition}

\noindent {\it Proof.} Our assumption implies that
\begin{equation}\left (\int_{G_{n,m}}|P_F(K)|^{-n}\,d\nu_{n,m}(F)\right )^{-\frac{1}{mn}}\ls
\left (\int_{G_{n,m}}|P_F(D)|^{-n}\,d\nu_{n,m}(F)\right )^{-\frac{1}{mn}}\end{equation}
By the linear invariance of $\Phi_{n-m}(D)$, for any $\tilde{D}=T(D)$ where $T\in GL(n)$ and $|\tilde{D}|=1$, we have
\begin{equation}\left (\int_{G_{n,m}}|P_F(D)|^{-n}\,d\nu_{n,m}(F)\right )^{-\frac{1}{mn}}=
|D|^{\frac{1}{n}}\left (\int_{G_{n,m}}|P_F(\tilde{D})|^{-n}\,d\nu_{n,m}(F)\right )^{-\frac{1}{mn}}.\end{equation}
Now, using H\"{o}lder's inequality we write
\begin{equation}\left (\int_{G_{n,m}}|P_F(\tilde{D})|^{-n}\,d\nu_{n,m}(F)\right )^{-\frac{1}{mn}}\ls
\left (\int_{G_{n,m}}|P_F(\tilde{D})|\,d\nu_{n,m}(F)\right )^{\frac{1}{m}}\end{equation}
From Aleksandrov's inequalites we have 
\begin{equation}\left (\int_{G_{n,m}}|P_F(\tilde{D})|\,d\nu_{n,m}(F)\right )^{\frac{1}{m}}\ls \omega_m^{\frac{1}{m}}w(\tilde{D})\ls c_2\sqrt{n/m}\frac{w(\tilde{D})}{\sqrt{n}}.\end{equation}
Taking into account Theorem \ref{th:main-3} we get
\begin{align}c\sqrt{n/m}|K|^{\frac{1}{n}}&\ls \left (\int_{G_{n,m}}|P_F(K)|^{-\frac{n}{m}}\,d\nu_{n,m}(F)\right )^{-\frac{1}{n}}
\ls \left (\int_{G_{n,m}}|P_F(D)|\,d\nu_{n,m}(F)\right )^{\frac{1}{m}}\\
\nonumber &= \left (\int_{G_{n,m}}|P_F(\tilde{D})|\,d\nu_{n,m}(F)\right )^{\frac{1}{m}}|D|^{\frac{1}{n}}
\ls c_2\sqrt{n/m}\frac{w(\tilde{D})}{\sqrt{n}}|D|^{\frac{1}{n}},\end{align}
and the result follows. \prend

\medskip

As a corollary we have:

\begin{theorem}\label{th:SP-2}Let $1\ls m\ls n-1$ and let $K$ and $D$ be two convex bodies in ${\mathbb R}^n$ such that
\begin{equation}|P_F(K)|\ls |P_F(D)|\end{equation} for every $F\in G_{n,m}$. Then,
\begin{equation}\label{eq:SP-2}|K|^{\frac{1}{n}}\ls c(\log n)\,|D|^{\frac{1}{n}},\end{equation}
where $c>0$ is an absolute constant.
\end{theorem}

\noindent {\it Proof.} If $D$ is in the minimal mean width position, we have (see \cite[Chapter 6]{AGA-book}) 
\begin{equation}w(\overline{D})\ls c_1\sqrt{n}\log n.\end{equation}
The result follows from Proposition \ref{prop:SP-1}. \prend

\section{Separation in the Busemann-Petty problem}\label{sect5}

For the proof of Theorem \ref{main-int} we need several definitions from convex geometry. A closed bounded set $K$ in $\R^n$ is called a star body if
every straight line passing through the origin crosses the boundary of $K$ at exactly two points different from the origin,
the origin is an interior point of $K$, and the Minkowski functional of $K$ defined by
\begin{equation}\label{eq:sep-2}\|x\|_K = \min\{a\gr 0:\ x\in aK\}\end{equation}
is a continuous function on $\R^n$.

The radial function of a star body $K$ is defined by
\begin{equation}\label{eq:sep-3}\rho_K(x) = \|x\|_K^{-1}, \qquad x\in \R^n,\ x\neq 0.\end{equation}
If $x\in S^{n-1}$ then $\rho_K(x)$ is the radius of $K$ in the direction of $x$.

We use the polar formula for volume of a star body
\begin{equation}\label{polar-volume}|K|=\int_{S^{n-1}} \|\theta\|_K^{-n} d\theta,
\end{equation}
where $d\theta$ stands for the uniform measure on the sphere with density 1.

The class of intersection bodies was introduced by Lutwak in \cite{Lutwak-1988}. Let $K, D$ be origin-symmetric star bodies in $\R^n.$ We say that $K$ is the
intersection body of $D$ and write $K=ID$ if the radius of $K$ in every direction is equal to the $(n-1)$-dimensional volume of the section of $L$ by the central
hyperplane orthogonal to this direction, i.e. for every $\xi\in S^{n-1},$
\begin{align}\label{eq:sep-4}
\rho_K(\xi) &= \|\xi\|_K^{-1} = |D\cap \xi^\bot|= \frac 1{n-1} \int_{S^{n-1}\cap \xi^\bot} \|\theta\|_D^{-n+1}d\theta \\
\nonumber &=\frac 1{n-1} R\left(\|\cdot\|_D^{-n+1}\right)(\xi),
\end{align}
where $R:C(S^{n-1})\to C(S^{n-1})$ is the {spherical Radon transform}
\begin{equation}\label{eq:sep-5}Rf(\xi)=\int_{S^{n-1}\cap \xi^\bot} f(x) dx,\qquad \hbox{for all}\; f\in C(S^{n-1}).\end{equation}
All bodies $K$ that appear as intersection bodies of different star bodies
form {the class of intersection bodies of star bodies}. A more general class of { intersection bodies}
is defined as follows. If $\mu$ is a finite Borel measure on $S^{n-1},$ then the spherical Radon transform
$R\mu$ of $\mu$ is defined as a functional on $C(S^{n-1})$ acting by
\begin{equation}\label{eq:sep-6}(R\mu, f)=(\mu, Rf)=\int_{S^{n-1}} Rf(x) d\mu(x),\qquad \hbox{for all}\; f\in C(S^{n-1}).\end{equation}
A star body $K$ in $\R^n$ is called an {\it intersection body} if $\|\cdot\|_K^{-1}=R\mu$
for some measure $\mu,$ as functionals on $C(S^{n-1}),$  i.e.
\begin{equation}\label{def-int}
\int_{S^{n-1}} \|x\|_K^{-1} f(x) dx = \int_{S^{n-1}} Rf(x)d\mu(x),\qquad \hbox{for all}\; f\in C(S^{n-1}).
\end{equation}
Intersection bodies played the key role in the solution of the Busemann-Petty problem.
\smallbreak
Recall that $d\sigma(x)=dx/|S^{n-1}|$ is the normalized uniform measure on the sphere, and denote by
$$M(K)=\int_{S^{n-1}} \|x\|_K d\sigma(x).$$
\bigbreak

\noindent {\bf Proof of Theorem \ref{main-int}.}   By \eqref{eq:sep-4}, the condition \eqref{sect1} can be written as
\begin{equation} \label{radon-buspetty}
R(\|\cdot\|_K^{-n+1})(\xi) \ls R(\|\cdot\|_L^{-n+1})(\xi) + (n-1)\e,
\ \hbox{for all}\; \xi\in S^{n-1}.
\end{equation}
Since $K$ is an intersection body, there exists a finite Borel measure
$\mu$ on $S^{n-1}$ such that $\|\cdot\|_K^{-1}= R\mu$ as functionals on $C(S^{n-1}).$
Together with \eqref{polar-volume}, \eqref{radon-buspetty} and the definition of $R\mu$, the latter implies that
\begin{align}\label{eq:sep-7} n|K| &= \int_{S^{n-1}} \|x\|_K^{-1}\|x\|_K^{-n+1}\ dx
= \int_{S^{n-1}} R\left(\|\cdot\|_K^{-n+1}\right)(\xi)\ d\mu(\xi)\\
\nonumber &\ls \int_{S^{n-1}} R\left(\|\cdot\|_L^{-n+1}\right)(\xi)\ d\mu(\xi) - (n-1)\e \int_{S^{n-1}}  d\mu(\xi)\\
\label{eq:sep-8} &= \int_{S^{n-1}} \|x\|_K^{-1}\|x\|_L^{-n+1}\ dx - (n-1)\e \int_{S^{n-1}} d\mu(x).
\end{align}
We estimate the first term in \eqref{eq:sep-8} using H\"older's inequality:
\begin{equation}\label{eq:sep-9}\int_{S^{n-1}} \|x\|_K^{-1}\|x\|_L^{-n+1}\ dx \ls  \left(\int_{S^{n-1}} \|x\|_K^{-n}\ dx\right)^{\frac1n}
\left( \int_{S^{n-1}} \|x\|_L^{-n}\ dx\right)^{\frac{n-1}n} = n |K|^{\frac1n} |L|^{\frac{n-1}n}.
\end{equation}
We now estimate the second term in \eqref{eq:sep-8} adding the Radon transform of  the unit constant
function under the integral ($R{\bf 1}(x)=\left|S^{n-2}\right|$ for every $x\in S^{n-1}$),
and using again the fact that $\|\cdot\|_K^{-1}=R\mu$:
\begin{align} \label{eq:sep-10}
(n-1)\e \int_{S^{n-1}} d\mu(x) &= \frac{(n-1)\e}{\left|S^{n-2}\right|} \int_{S^{n-1}} R1(x)\ d\mu(x)
=\frac{(n-1)\e}{\left| S^{n-2} \right| } \int_{S^{n-1}} \|x\|_K^{-1}\ dx\\
\nonumber &\gr c_1\e \frac{(n-1)|S^{n-1}|}{\left| S^{n-2} \right|} \frac{1}{M(\overline{K})}|K|^{\frac 1n}\\
\label{eq:sep-11} &\gr c_2\e \sqrt{n}  \frac{1}{M(\overline{K})} |K|^{\frac 1n},
\end{align}
since
\begin{equation}\int_{S^{n-1}}\|x\|_K^{-1}\,d\sigma(x)\ge \left (\int_{S^{n-1}}\|x\|_Kd\sigma(x)\right )^{-1}=\frac{1}{M(\overline{K})}|K|^{\frac{1}{n}},\end{equation}
by Jensen's inequality, homogeneity, and
\begin{equation}\label{eq:sep-12}\left|S^{n-2}\right|= \frac{2\pi^{\frac{n-1}2}}{\Gamma(\frac{n-1}2)} \qquad {\rm and}\qquad
\left|S^{n-1}\right|= \frac{2\pi^{\frac{n}2}}{\Gamma(\frac{n}2)}.\end{equation}
Combining \eqref{eq:sep-11} with \eqref{eq:sep-8} and  \eqref{eq:sep-9}, we get
\begin{equation}\label{eq:sep-13}n |K| \ls n |K|^{\frac1n} |L|^{\frac{n-1}n} - c_2\e \sqrt{n}  \frac{1}{M(\overline{K})} |K|^{\frac 1n}
\end{equation}
and, after dividing by $n|K|^{1/n},$ the proof is complete. \prend

\bigbreak

Separation implies a volume difference inequality.

\begin{corollary}\label{cor:sep-2}Let $L$ be any origin-symmetric star body in $\R^n,$ and let $K$ be an origin-symmetric
intersection body, which is a dilate of an isotropic body. Suppose that
$$\min_{\xi\in S^{n-1}} \left(|L\cap \xi^\bot|-|K\cap \xi^\bot|\right) > 0.$$
Then
$$|L|^{\frac{n-1}n} - |K|^{\frac{n-1}n} \gr c_2 \frac{1}{\sqrt{n}M(\overline{K})}
\min_{\xi\in S^{n-1}} \left(|L\cap \xi^\bot|-|K\cap \xi^\bot|\right).$$
\end{corollary}

\begin{remark}\label{rem:sep-3}\rm It was proved in \cite{GM} that there exists a constant $c>0$ such that
for any $n\in \N$ and any origin-symmetric isotropic convex body $K$ in $\R^n$
\begin{equation}\label{eq:sep-1}\frac{1}{M(K)} \gr\ c_1 \frac{n^{1/10}L_K}{\log^{2/5}(e+n)}\gr c_2\frac{n^{1/10}}{\log^{2/5}(e+n)}.\end{equation}
Also, if $K$ is convex, has volume $1$ and is in the minimal mean width position then we have
\begin{equation}\label{eq:sep-22}\frac{1}{M(K)} \gr\ c_3\frac{\sqrt{n}}{\log (e+n)}.\end{equation}
Inserting these estimates into Theorem \ref{main-int} and Corollary \ref{cor:sep-2} we obtain estimates independent
from the bodies.
\end{remark}

\bigbreak

\bigskip

\bigskip

\noindent {\bf Acknowledgements.} We would like to thank Silouanos Brazitikos for useful suggestions that helped us to simplify
the proof and improve the estimate in Theorem \ref{th:intro-S}. The second named author
was partially supported by the US National Science Foundation grant DMS-1265155.

\bigskip

\bigskip

\footnotesize
\bibliographystyle{amsplain}

\begin{thebibliography}{100}
\footnotesize

\bibitem{AGA-book}
\textrm{S.\ Artstein-Avidan, A.\ Giannopoulos and V.\ D.\ Milman},
{\sl Asymptotic Geometric Analysis, Vol. I}, Mathematical Surveys and Monographs {\bf 202}, Amer. Math. Society (2015).
\bibitem{Ball-1989} {\rm K.\ M.\ Ball}, {\sl Volumes of sections of cubes and related problems}, in Geom. Aspects of Funct.
Analysis, Lecture Notes in Mathematics {\bf 1376}, Springer, Berlin (1989), 251-260.
\bibitem{Ball-1991} {\rm K.\ M.\ Ball}, {\sl Shadows of convex bodies}, Trans.  Amer. Math. Soc. {\bf 327} (1991), 891-901.
\bibitem{Ball-1991c} {\rm K.\ M.\ Ball}, {\sl Volume ratios and a reverse isoperimetric inequality}, J. London Math. Soc. (2)
{\bf 44} (1991), 351-359.
\bibitem{Bourgain-1991} {\rm J.\ Bourgain}, {\sl On the distribution of polynomials on high dimensional convex sets}, in Geom. Aspects of Funct.
Analysis, Lecture Notes in Mathematics {\bf 1469}, Springer, Berlin (1991), 127-137.
\bibitem{BZ} {\rm J.~Bourgain and Gaoyong Zhang}, {\sl On a generalization of the
Busemann-Petty problem}, Convex geometric analysis (Berkeley, CA, 1996),
65--76, Math. Sci. Res. Inst. Publ., 34, Cambridge Univ. Press,
Cambridge, 1999.
\bibitem{BGVV-book} \textrm{S.\ Brazitikos, A.\ Giannopoulos, P.\ Valettas and B-H.\ Vritsiou},
{\sl Geometry of isotropic convex bodies}, Mathematical Surveys and Monographs {\bf 196}, Amer. Math. Society (2014).
\bibitem{Busemann-Petty-1956} {\rm H.\ Busemann and C.\ M.\ Petty},
{\sl Problems on convex bodies}, Math. Scand. {\bf 4} (1956), 88-94.
\bibitem{Burago-Zalgaller-book} {\rm Y.\ D.\ Burago and V.\ A.\ Zalgaller}, {\sl Geometric Inequalities}, Springer Series in
Soviet Mathematics, Springer-Verlag, Berlin-New York (1988).
\bibitem{Chasapis-Giannopoulos-Liakopoulos-2015} {\rm G.\ Chasapis, A.\ Giannopoulos and D-M.\ Liakopoulos},
{\sl Estimates for measures of lower dimensional sections of convex bodies}, Preprint.
\bibitem{Dafnis-Paouris-2012} {\rm N.\ Dafnis and G.\ Paouris}, {\sl Estimates for the affine and dual
affine quermassintegrals of convex bodies}, Illinois J.\ of Math.\ {\bf 56} (2012), 1005-1021.
\bibitem{Gardner-book}{\rm R.\ J.\ Gardner}, {\sl Geometric Tomography}, Second Edition
Encyclopedia of Mathematics and its Applications {\bf 58}, Cambridge
University Press, Cambridge (2006).
\bibitem{GM} {\rm A.\ Giannopoulos and E.\ Milman}, {\sl M-estimates for isotropic convex bodies and their Lq-centroid bodies},
Geometric aspects of functional analysis, 159--182, Lecture Notes in Math., 2116, Springer, Cham, 2014.
\bibitem{GZ}{\rm P. Goodey and Gaoyong Zhang,} {\sl Inequalities between
projection functions of convex bodies}, Amer. J. Math. {\bf 120} (1998),
345-367.
\bibitem{Grinberg-1990} {\rm E.\ L.\ Grinberg},  {\sl Isoperimetric
inequalities and identities for $k$-dimensional cross-sections of convex bodies}, Math. Ann. {\bf 291} (1991), 75-86.
\bibitem{Klartag-2006} \textrm{B.\ Klartag}, {\sl On convex perturbations with a bounded isotropic constant},
Geom.\ Funct.\ Anal.\ {\bf 16} (2006), 1274-1290.
\bibitem{Klartag-EMilman-2012} \textrm{B.\ Klartag and E.\ Milman},
{\sl Centroid Bodies and the Logarithmic Laplace Transform -- A Unified Approach},
J.\ Funct.\ Anal.\ \textbf{262} (2012), 10--34.
\bibitem{Koldobsky-book} {\rm A.\ Koldobsky}, {\sl Fourier analysis in convex geometry},
Mathematical Surveys and Monographs {\bf 116}, Amer. Math. Society (2005).
\bibitem{K??} {\rm A.\ Koldobsky}, {\sl Stability in the Busemann-Petty and Shephard problems},
Adv. Math. {\bf 228} (2011), 2145--2161.
\bibitem{Koldobsky-3} {\rm A.~Koldobsky}, {\sl Stability and separation in volume comparison problems,}
Math. Model. Nat. Phenom. {\bf 8} (2012), 159--169.
\bibitem{Koldobsky-4} {\rm A.~Koldobsky}, {\sl Isomorphic Busemann-Petty problem for sections of proportional
dimensions,} Adv. in Appl. Math.  {\bf 71} (2015), 138--145.
\bibitem{Lutwak-1984}{\rm E. Lutwak}, {\sl A general isepiphanic
inequality}, Proc. Amer. Math. Soc. {\bf 90} (1984), 415--421.
\bibitem{Lutwak-1988b}{\rm E. Lutwak}, {\sl Inequalities for Hadwiger's
harmonic Quermassintegrals}, Math. Annalen {\bf 280} (1988),
165--175.
\bibitem{Lutwak-1988} {\rm E. Lutwak}, {\sl Intersection bodies and dual
mixed volumes}, Advances in Math. {\bf 71} (1988), 232-261.
\bibitem{Paouris-Pivovarov-2013} {\rm  G.\ Paouris and P.\ Pivovarov},
{\sl Small-ball probabilities for the volume of random convex sets}, Discrete Comput. Geom. {\bf 49} (2013), 601-646.
\bibitem{Petty-1967} {\rm C.\ M.\ Petty}, {\sl Projection bodies}, Proc. Colloquium
Convexity, Copenhagen 1965, Kobenhavns Univ. Mat. Inst., 1967, 234-241.
\bibitem{Pisier-1989} {\rm G.\ Pisier}, {\sl A new approach to several results of V. Milman},
J. Reine Angew. Math. {\bf 393} (1989), 115-131.
\bibitem{Schneider-1967} {\rm R.\ Schneider}, {\sl Zu einem Problem von Shephard
\"{u}ber die Projektionen konvexer K\"{o}rper}, Math. Z. {\bf 101}
(1967), 71-82.


\bibitem{Schneider-book} {\rm R.\ Schneider}, {\sl Convex Bodies: The Brunn-Minkowski Theory},
Second expanded edition. Encyclopedia of Mathematics and Its Applications 151, Cambridge University Press, Cambridge, 2014.
\bibitem{Shephard-1964} {\rm G.\ C.\ Shephard}, {\sl Shadow systems of convex bodies},
Israel J. Math. {\bf 2} (1964), 229-236.
\end{thebibliography}

\bigskip

\bigskip

\thanks{\noindent {\bf Keywords:}  Convex bodies; Busemann-Petty problem; Shephard problem;
Affine and dual affine quermassintegrals; Intersection body; Isotropic convex body.}

\smallskip

\thanks{\noindent {\bf 2010 MSC:} Primary 52A20; Secondary 46B06, 52A23, 52A40.}

\bigskip

\bigskip

\noindent \textsc{Apostolos \ Giannopoulos}: Department of
Mathematics, University of Athens, Panepistimioupolis 157-84,
Athens, Greece.

\smallskip

\noindent \textit{E-mail:} \texttt{apgiannop@math.uoa.gr}

\bigskip

\noindent \textsc{Alexander \ Koldobsky}: Department of
Mathematics, University of Missouri, Columbia, MO 65211.

\smallskip

\noindent \textit{E-mail:} \texttt{koldobskiya@missouri.edu}

\bigskip

\end{document}